\documentclass[letterpaper, 10 pt, conference]{ieeeconf}

\usepackage{graphicx}
\usepackage{amsmath}
\usepackage{amssymb}
\usepackage{CJK}
\usepackage{color}
\usepackage{dsfont}
\usepackage{subfigure}
\usepackage{booktabs}
\usepackage{microtype}
\usepackage{mathrsfs}
\usepackage{bm}
\usepackage{cite}
\usepackage{url}
\usepackage[colorlinks=true,allcolors=blue]{hyperref}

\DeclareMathOperator{\Tr}{Tr}
\IEEEoverridecommandlockouts                              
\title{A Performance and Stability Analysis of Low-inertia Power Grids with Stochastic System Inertia
}

\author{Yi Guo, Tyler H. Summers
\thanks{Y. Guo and T.H. Summers are with the Department
of Mechanical Engineering, The University of Texas at Dallas, Richardson,
TX, USA, email: \{yi.guo2,tyler.summers\}@utdallas.edu.}
\thanks{This material is based on work supported by the National Science Foundation under grants CNS-1566127 and CMMI-1728605.}
\thanks{This preprint is an extended version of a conference paper \cite{guo2019performance} accepted to 2019 American Control Conference, Philadelphia, USA, July 10-12, 2019.}
}

\begin{document}
\maketitle
\thispagestyle{empty}
\pagestyle{empty}
\begin{abstract}
Traditional synchronous generators with rotational inertia are being replaced by low-inertia renewable energy resources (RESs) in many power grids and operational scenarios. Due to emerging market mechanisms, inherent variability of RESs, and existing control schemes, the resulting system inertia levels can not only be low but also markedly time-varying. In this paper, we investigate performance and stability of low-inertia power systems with stochastic system inertia. In particular, we consider system dynamics modeled by a linearized stochastic swing equation, where stochastic system inertia is regarded as multiplicative noise. The $\mathcal{H}_2$ norm is used to quantify the performance of the system in the presence of persistent disturbances or transient faults.
The performance metric can be computed by solving a generalized Lyapunov equation, which has fundamentally different characteristics from systems with only additive noise.
For grids with uniform inertia and damping parameters, we derive closed-form expressions for the $\mathcal{H}_2$ norm of the proposed stochastic swing equation. The analysis gives insights into how the $\mathcal{H}_2$ norm of the stochastic swing equation depends on 1) network topology; 2) system parameters; and 3) distribution parameters of disturbances. A mean-square stability condition is also derived. Numerical results provide additional insights for performance and stability of the stochastic swing equation.

\end{abstract}

\color{black}
\section{Introduction}
Environmental and sustainability concerns are forcing unprecedented charges in the modern electric power system. The continued replacement of traditional synchronous generators by renewable energy sources (RESs) in power systems is raising concerns about their stability. As the penetration levels of RESs reach substantial fractions of total generation, power systems will require more low-inertia RESs to participate in frequency and voltage control. The inherent variability of RESs can produce high amplitude and persistent disturbances, which may adversely affect stability. Due to emerging market mechanisms and deregulated existing control schemes of RESs, the resulting system inertia levels can not only be low but also markedly time-varying. This system-level inertia variation together with unpredictable RESs and net loads make frequency control and power system stabilization more challenging. Future power systems need more sophisticated stochastic dynamic models and stochastic control methods for stability and performance analysis.

Many issues of power system stability have been well studied in recent decades through mathematical analysis and computational techniques \cite{chiang1987foundations,chiang2011direct,kundur2004definition,bergen1981structure,tegling2015price,pirani2017system,motter2013spontaneous,dorfler2012synchronization,luo2018coordinated,ahmadyar2018framework}. However, many of the underlying assumptions and models must be called into question in the context of low- and variable-inertia power systems. The decreasing system inertia results in higher rate of change of frequency, which requires controllers to respond faster to the system dynamics. In addition, the intermittent disturbances from RESs are spatially distributed over power networks, which requires that performance and stability analysis consider various grid topologies and the diverse dynamics of various grid-connected components.

Many challenges and related solutions for low-inertia power grid stability have been highlighted and discussed in \cite{ulbig2013impact,ulbig2015analyzing,kroposki2017achieving,winter2015pushing,tielens2012grid}. Recently, stability analysis and control techniques have been proposed and have demonstrated their effectiveness for system stabilization \cite{beck2007virtual,d2013virtual,zhong2012control,dorfler2016breaking,bevrani2014virtual,tegling2015price,pirani2017system,dorfler2012synchronization,zhu2018stability,song2017network,colombino2017global,torres2015synchronization,sinha2017uncovering,zhao2014design,mallada2017optimal,curi2017control,li2016connecting,ahmadyar2018framework,mevsanovic2016comparison,ahmadyar2016assessment}. Virtual inertia emulation is an approach to control the terminal behavior of 
inverter-interfaced RESs to mimic inertial response of conventional synchronous machines \cite{beck2007virtual,d2013virtual,zhong2012control,dorfler2016breaking,bevrani2014virtual}. Recent works on virtual oscillators have shown that oscillator-based control strategies 
have advantages in faster response and global convergence \cite{colombino2017global,torres2015synchronization,sinha2017uncovering}. Also, distributed control methodologies have been implemented on generation and load sides to provide effective ways for frequency stabilization \cite{dorfler2016breaking,zhao2014design,mallada2017optimal,curi2017control,li2016connecting}. Recent stability analyses have demonstrated that the efforts to maintain synchronous stability in low-inertia power grids depend on grid structure, node dynamics and coupling strength \cite{tegling2015price,pirani2017system,dorfler2012synchronization,zhu2018stability,song2017network}. 

Overall, this line of research has explored useful solutions for the well recognized issues in low-inertia power systems. However, a large portion of unpredictable RESs will possible cause the system inertia to be not just low but also significantly time-varying. A few recent analysis consider the effects of system inertia variability raised by RESs \cite{ahmadyar2018framework,ahmadyar2016assessment,mevsanovic2016comparison,poolla2017optimal}, but none of works explicitly consider a rigorous performance and stability analysis with respect to a stochastic system inertia model.  As the inverter-based RESs dominated the generation, some RESs are required to 
participate in frequency and voltage control in a low/no-inertia power grid. The time-varying system inertia profiles and the heterogeneous allocation of inertia can lead to destabilizing effects, which complicates and challenges stability analysis and stabilization control in power systems. 

In this paper, we investigate performance and stability of low-inertia power systems with stochastic system inertia. In particular, we consider system dynamics modeled by a linearized stochastic swing equation, where stochastic system inertia is regarded as multiplicative noise. The main contributions are as follows:
\begin{itemize}
\item We consider the frequency dynamics of a low-inertia power grid as a stochastic linear system with both multiplicative (due to inertia variations) and additive noise (due to power injection disturbances), which connects stability analysis of a linearized swing equation with a generalized Lyapunov equation.

\item We quantify the system $\mathcal{H}_2$ norm in terms of various outputs, which measures system performance in the presence of multiplicative and additive disturbances. The closed-form $\mathcal{H}_2$ norm of a homogeneous power grid is derived and discussed. In contrast to existing work with only additive disturbances, we observe that the system $\mathcal{H}_2$ norm depends in fundamentally different ways on 1) network topology; 2) system parameters; 3) distribution parameters of disturbances, and is sensitive to system inertia variation for certain outputs. 

\item A mean-square stability condition is also derived for the stochastic linear system, which indicates that a low-inertia grid can be destabilized in a second-order sense by inertia variability. Numerical results also indicate that a lower inertia grid with larger system inertia variance is less robust to disturbances. 

\end{itemize}

The rest of the paper is organized as follows: Section II introduces a stochastic linearized swing equation to model power systems with stochastic system inertia. Section III develops the system $\mathcal{H}_2$ norm of the stochastic swing equation for three particular outputs and derives a second-order stability condition. Section IV presents some numerical results to illustrate the theory. Section V concludes the paper and points out several possible future research directions.



\section{Problem formulation}

\subsection{System modelling}
Consider a power network $\mathcal{G}$ with $N+1$ nodes (buses) $\mathcal{V} = \{0\} \cup \{1, \ldots, N\}$ connected by a set of edges (transmission lines) $\mathcal{E} \subset \{\mathcal{V} \times \mathcal{V}\}$. 
We assume the power network is Kron-reduced \cite{dorfler2013kron}, where each node represents an equivalent generator\footnote{The equivalent generator can be interpreted as the aggregate terminal interaction dynamics of a power system sub-area. Note that the equivalent generator here characterizes the dynamics of a group of grid-connected components (e.g., traditional generators and inverter-based generators) as a synchronous machine.} with state variables (e.g., voltage magnitude $|V_i|$ and voltage angle $\theta_i$) and parameters (e.g., rotational inertia $M_i$, damping coefficient $\beta_i$). The index $0$ is reserved for the grounded node/bus in a Kron-reduced power network. We consider the following swing equation to model the generator dynamics of each bus $i$

\begin{equation}\label{ClassSwingEquation}
M_i\ddot{\theta}_i + \beta_i\dot{\theta}_i = P_{m,i} - P_{e,i}, \forall i = 0,1, \ldots , N,
\end{equation}
where $P_{m,i}$ refers to the mechanical power from the synchronous generator, and $P_{e,i}$ represents the electrical power injection of the generator. The nominal power injection at bus $i$ is given by the power flow equations
\begin{equation}\label{nonlinearPowerFlow}
\begin{aligned}
P_{e,i} = \overline{g}_{ii}|V_i|^2 + \sum_{(i,j)\in \mathcal{E}}g_{ij}|V_i||V_j|\cos(\theta_i-\theta_j) ~~~~~~~~~ \\
+ \sum_{(i,j) \in \mathcal{E}} b_{ij}|V_i||V_j|\sin(\theta_i - \theta_j).
\end{aligned}
\end{equation}
where $g_{ij}$ and $b_{ij}$ denote the line conductance and the line susceptance, respectively. The shunt capacity at bus $i$ is $\overline{g}_{ii}$. The admittance matrix $\mathbf{Y} \in \mathbf{C}^{(N+1) \times (N+1)}$ has elements
\begin{equation}\label{admittanceMatrix}\nonumber
Y_{ij} = \left\{ \begin{array}{ll}
\sum_{l \sim i} (g_{il} - \mathbf{j} b_{il}) + \overline{g}_{ii} & \textrm{if $i=j$}\\
-(g_{ij} - \mathbf{j}b_{ij}) & (i,j) \in \mathcal{E} \\
0 &  (i,j) \notin \mathcal{E}.
\end{array} \right.
\end{equation}
The admittance matrix can be written in compact matrix form
\begin{equation}\nonumber
\mathbf{Y} = (L_G + \overline{g}) - \mathbf{j}L_B,
\end{equation}
where $L_G$ and $L_B$ are the conductance matrix and the susceptance matrix, and $\overline{g}:= \textrm{diag}\{\overline{g}_{ii}\}$ is a diagonal matrix of bus shunt capacitors. The Laplacian matrices $L_B$ and $L_G$ comprise the weights of the line susceptance $b_{ij}$ and the line conductance $g_{ij}$ in the Kron-reduced network, respectively. 

The dynamic model \eqref{ClassSwingEquation}-\eqref{nonlinearPowerFlow} is often linearized around an operating point, which allows the study of system response in the presence of small faults or persistence disturbances around the linearization point. We use the linearized ``DC power flow model'', neglecting the line reactance, to approximate \eqref{nonlinearPowerFlow}, which assumes $|V_i| \approx 1$ and $|\theta_i - \theta_j| \ll 1$. Detailed analysis and applications of the DC power flow approximation are discussed in \cite{purchala2005usefulness}. Then the system dynamics \eqref{ClassSwingEquation} becomes
\begin{equation}\label{LinearizedSwingEquation}
M_i\ddot{\theta_i} + \beta_i\dot{\theta_i} \approx - \sum_{(i,j)\in \mathcal{E}} b_{ij}(\theta_i - \theta_j) + P_{m,i}.
\end{equation}
We then shift the equilibrium point of \eqref{LinearizedSwingEquation} to the origin and write it in state-space form
\begin{equation}\label{LTI ODE expression}
\begin{aligned}
\dot{\theta} & = \omega \\
M\dot{\omega} & = - L_B\theta -D\omega + W,
\end{aligned}
\end{equation}
where $\theta = [\theta_1,\ldots,\theta_N]^\intercal$ and $\omega = [\omega_1,\ldots,\omega_N]^\intercal$. The inertia matrix and damping matrix are defined as $M := \textrm{diag}\{M_i\}$, $D := \textrm{diag}\{\beta_i\}$. The standard approach to analyze  \eqref{LinearizedSwingEquation} considers various disturbances as additive noise $W$ driving the system away from its current equilibrium point. The remainder of this paper considers both \emph{multiplicative} and \emph{additive} disturbances in \eqref{LTI ODE expression}.

\subsection{Frequency dynamics with multiplicative and additive noise}
Here we consider inertia variations caused by RESs, which are modeled by treating the system inerta matrix as multiplicative noise rather than simply a constant. The inertia parameter at each node $M_i$ can be modeled as a independent Wiener processes on a probability space $(\Omega, \mathcal{F},\mu)$ with mean $\bar{M}_i$ and variance $\sigma^2_i$ \cite{oksendal2003stochastic}. The nominal inertia matrix is defined as $\mathcal{M} := \textrm{diag}\{\bar{M}_i\}$, 
and we rewrite \eqref{LTI ODE expression} as a stochastic linear system with both additive and multiplicative noise
\begin{equation}\label{LTVstatespace}
\scalebox{0.85}[1]{$
\begin{bmatrix}
\dot{\theta}\\
\dot{\omega}
\end{bmatrix}=
\begin{bmatrix}
0 & I\\
-(\mathcal{M}^{-1} + \bm{\delta} \mathcal{M}^{-1})L_B & -(\mathcal{M}^{-1} +  \bm{\delta} \mathcal{M}^{-1})D
\end{bmatrix}
\begin{bmatrix}
\theta \\
\omega
\end{bmatrix}+
\begin{bmatrix}
0\\
\eta I
\end{bmatrix}W.
$}
\end{equation}
The matrix $\mathcal{M}^{-1} := \textrm{diag}\{\hat{M}^{-1}_i\}$ collects the mean values of the inverse distribution of $\bar{M}_i$ on the diagonal, and the additive noise $W$ represents independent white-noise with zero mean and unit variance, scaled by $\eta$. Each diagonal element of the matrix $\bm{\delta}\mathcal{M}^{-1}: = \textrm{diag}\{ \delta_i \}$ is modeled as an independent Wiener process with zero mean and variance $\hat{\sigma}^2_i$. For simplicity, we write \eqref{LTVstatespace} in a generalized form with outputs as a multi-input multi-output stochastic linear system with multiplicative noise
\begin{equation}\label{LTVswingequation}
\begin{aligned}
& \dot{x} = A_0x + \sum_{i = 1}^{N+1}  A_i \delta_i x + BW,\\
& y = Cx,
\end{aligned}
\end{equation}
where,
\begin{equation*} A_0 = 
\begin{bmatrix}
0 & I\\
-\mathcal{M}^{-1}L_B & -\mathcal{M}^{-1}D 
\end{bmatrix},
A_i =
\begin{bmatrix}
0 & 0\\
-R_iL_B & -R_iD
\end{bmatrix},
\end{equation*}
\begin{equation*}
B = 
\begin{bmatrix}
0\\
\eta I
\end{bmatrix},
x = 
\begin{bmatrix}
\theta^\intercal, \omega^\intercal
\end{bmatrix}^\intercal.
\end{equation*}
The matrix $A_0 \in \mathbf{R}^{2(N+1)}$ characterizes the nominal system with average inertia and damping ratio. We define an \emph{inertia disturbance allocation} matrix $R_i \in \mathbf{R}^{N+1}$ in $A_i \in \mathbf{R}^{2(N+1)}$ associated with each bus $i$. The elements in $R_i$ are all zeros except for one diagonal element $r_{ii} = 1$, which maps the corresponding inertia disturbance $\delta_i$ onto bus $i$. If the inertia variation at bus $i$ is insignificant, we set $R_i = 0$ to remove the inertia disturbance at $i$th bus. In the following, we refer to the stochastic system input/output mapping \eqref{LTVswingequation} as $\Sigma = (A_0, A_i, B, C)$. This stochastic linear dynamical model allows us to investigate the effects of phase angle deviations and frequency changes in the presence of both additive and multiplicative disturbances around the original operating point. To assess the system stability and evaluate the performance of \eqref{LTVswingequation}, we consider the following three outputs \cite{pirani2017system}:

\subsubsection{Phase cohesiveness} This output quantifies real power losses due to phase differences caused by the fluctuations from the nominal operating points \cite{tegling2015price}. The resistive losses on transmission lines during transients or due to persistent disturbances can be expressed in terms of Laplacian matrix $L_G$
\begin{equation}\nonumber
P_{\textrm{loss}} = \sum_{(i,j) \in \mathcal{E}}g_{ij}(\theta_i - \theta_j)^2 = \theta^\intercal L_G \theta. 
\end{equation}
Expressed in terms of the output, we have $P_{\textrm{loss}} = y^\intercal y$ with
\begin{equation*}
y = 
\begin{bmatrix}
L_G^{\frac{1}{2}} & 0
\end{bmatrix}
x.
\end{equation*}
It is worth to note that the stochastic linear system \eqref{LTVswingequation} driven from the linearized swing equation neglects the line resistances. The output matrix with conductance matrix $L_G$ can capture resistive losses arising from \eqref{LTVswingequation}.

\subsubsection{Frequency} To quantify frequency deviations due to faults or disturbances, we have the output $y =  [0, I]x$.

\subsubsection{Phase cohesiveness \& frequency} This output quantifies both phase and frequency performance with the output matrix
\begin{equation*}
y =
\textrm{diag}\{L_G^{\frac{1}{2}},\kappa I\}x,
\end{equation*}
where $\kappa \in \mathbf{R}_{+}$ trades off phase angle and frequency deviations.

In this paper, we use the system $\mathcal{H}_2$ norm to quantify the system performance under above three outputs. The system $\mathcal{H}_2$ norm is the root-mean-square value of the output when the system is driven by multiplicative and additive noise inputs.  It has been widely studied in power system models with additive noise \cite{pirani2017system,poolla2017optimal,tegling2015price,jiang2017performance}.

\section{Coherency performance metric}
\subsection{System reduction}
The stochastic swing equation model \eqref{LTVswingequation} and associated performance outputs consist of two key Laplacian matrices: the conductance matrix $L_G$ and the susceptance matrix $L_B$, which both have a zero eigenvalue associated with the eigenvector $\bm{1}$. The Laplacian structure implies \eqref{LTVswingequation} is not asymptotically stable, but the subspace corresponding to the zero eigenvalue does not appear in the output \cite{tegling2015price}. The physical interpretation of this zero eigenvalue is that we lack a grounded reference bus in the power network. Therefore, we consider the grounded Laplacians by deleting the $k$th row and column of $L_G$ and $L_B$, respectively, yielding $\tilde{L}_G$ and $\tilde{L}_B$. The states in the reduced system $\tilde{\theta}$ and $\tilde{\omega}$ are obtained by simply removing the $k$th element of the original vectors, which can be interpreted as grounding bus $k$, with $\theta_k$ and $\omega_k$ forced to zero for the voltage references. Then the grounded system $\tilde{\Sigma}$ can be expressed as
\begin{equation}\label{LTVstateform_reduced}
\begin{aligned}
& \dot{\tilde{x}} = \tilde{A}_0\tilde{x} + \sum_{i = 1}^N  \tilde{A}_i \delta_i\tilde{x} + \tilde{B}\tilde{W},\\
& \tilde{y} = \tilde{C}\tilde{x},
\end{aligned}
\end{equation}
where,
\begin{equation*} \tilde{A}_0 = 
\begin{bmatrix}
0 & I\\
-\mathcal{\tilde{M}}^{-1}\tilde{L}_B & -\mathcal{\tilde{M}}^{-1}\tilde{D} 
\end{bmatrix},
\tilde{A}_i =
\begin{bmatrix}
0 & 0\\
-\tilde{R}_i\tilde{L}_B & -\tilde{R}_i\tilde{D}
\end{bmatrix},
\end{equation*}
\begin{equation*}
\tilde{B} = 
\begin{bmatrix}
0\\
\eta I
\end{bmatrix},
\tilde{x} = 
\begin{bmatrix}
\tilde{\theta}^\intercal, \tilde{\omega}^\intercal
\end{bmatrix}^\intercal.
\end{equation*}
    The $k$th row and column of matrices $\mathcal{M}$, $R_i$ and $D$ are also discarded in the reduced power system \eqref{LTVstateform_reduced}, which are then written as $\tilde{\mathcal{M}}$, $\tilde{R}_i$ and $\tilde{D}$, respectively. We assume the power network $\mathcal{G}$ in our problem is connected so that the grounded Laplacians $\tilde{L}_G$ and $\tilde{L}_B$ are symmetric positive definite. Thus, all eigenvalues of matrix $\tilde{A}_0$ are located in the open left half of the complex plane. We will detail the system stability with finite $\mathcal{H}_2$ norm in the present of multiplicative noise in the next subsection. 

\subsection{Performance metric and stability conditions}
The system $\mathcal{H}_2$ norm is the root-mean-square value of the output when the system is driven by multiplicative and additive noise inputs.
For the system $\tilde{\Sigma} = (\tilde{A}_0, \tilde{A}_i, \tilde{B}, \tilde{C})$, the squared $\mathcal{H}_2$ norm is given by
\begin{equation}\label{H2observability}
\|\tilde{\Sigma}\|_{\mathcal{H}_2}^2 = \Tr(\tilde{B}^\intercal\tilde{Q}\tilde{B}),
\end{equation}
where $\tilde{Q}$ is the observability Gramian \cite{damm2004rational}. The observability Gramian $\tilde{Q}$ can be interpreted as the steady-state output covariance \cite{verriest2008time}, which has a stochastic integral expression \cite{zhang2002h2}. When a finite positive semidefinite observability Gramian $\tilde{Q}$ exists, it can be attained by solving the following generalized Lyapunov equation \eqref{generalizedLyapunov} \cite{kleinman1969stability,damm2004rational,zhang2002h2}
\begin{equation}\label{generalizedLyapunov}
\tilde{A}_0^\intercal\tilde{Q} + \tilde{Q}\tilde{A}_0 + \sum_{i=1}^N \hat{\sigma}_i^2\tilde{A}_i^\intercal\tilde{Q}\tilde{A}_i = - \tilde{C}^\intercal \tilde{C},
\end{equation}
which also implies the system is second-moment bounded (i.e., mean-square stable), resulting in the finite $\mathcal{H}_2$ norm.

The Gramian obtained from this generalized Lyapunov equation is used to compute the $\mathcal{H}_2$ norm  of system \eqref{LTVstateform_reduced}, which explicitly incorporates both multiplicative and additive noise. The multiplicative noise can be removed by letting $\hat{\sigma}_i^2 = 0$ in \eqref{generalizedLyapunov}, yielding a standard Lyapunov equation. 
In contrast to systems with only additive noise, there are differing notions of stability when multiplicative noise is present. In particular, even when the mean value of the state is stable, (i.e., $\tilde{A}_0$ is stable), the covariance of the state may be unstable due to the multiplicative noise, in which case the $\mathcal{H}_2$ norm becomes infinite. When the multiplicative noise variances are sufficiently small, the system will have second-moment bounded (i.e., be mean-square stable). We have the following second-order stability definition, which is equivalent to existence of a finite positive semidefinite solution to the generalized Lyapunov equation and finiteness of the corresponding $\mathcal{H}_2$ norm \cite{damm2004rational}.

\textbf{Definition 1: [Second-moment boundedness].} The system \eqref{LTVstateform_reduced} is called second moment bounded, or mean square stable, if there exists a positive constant $\alpha$ such that
\begin{equation*}
    \lim_{t \to \infty} \mathds{E}[\tilde{x}(t)^\intercal \tilde{x}(t)] \leq \alpha,~~\forall \tilde{x}(0) \in \mathbf{R}^N. 
\end{equation*}


\color{black}
\textbf{Remark 1:} \textbf{[Computation of \bm{$\mathcal{H}_2$} norm for non-homogeneous power grids].} The generalized Lyapunov equation \eqref{generalizedLyapunov} is linear in $\tilde{Q}$ and can be solved directly using vectorization and Kronecker products, yielding
\begin{equation}\nonumber
\scalebox{0.93}[1]{$
\textrm{vec}(\tilde{Q}) = - \left( \tilde{A}_0^\intercal \otimes I + I \otimes \tilde{A}_0^\intercal + \sum_{i=1}^N \tilde{A}_i^\intercal \otimes \tilde{A}_i^\intercal \right)^{-1}\textrm{vec}(\tilde{C}^\intercal \tilde{C}).
$}
\end{equation}
For the standard Lyapunov equation, factorization methods can be used to exploit the structure of the equation and achieve superior computational complexity. However, these methods cannot be easily applied to the generalized Lyapunov equation, and alternative methods have been studied \cite{benner2011lyapunov,damm2008direct,bai2006projection,boyd1994linear}, e.g., using Krylov subspaces, semidefinite programming, or differential equations for the state covariance matrix.

To gain additional insights into the effects of multiplicative noise on power networks, we now consider computation of the  $\mathcal{H}_2$ norm for power networks with homogeneous nominal inertia and damping. The inverse inertia perturbation at each bus is an independent stochastic process, which has identical mean and variance $\hat{M}_i = \hat{M}$, $\hat{\sigma}_i^2 = \hat{\sigma}^2, \forall i$. The damping ratio is also assumed to be identical $D = \beta I$. We will derive a closed-form expression for the $\mathcal{H}_2$ norm of the stochastic system \eqref{LTVstateform_reduced}, which allows us to highlight several insights regarding system performance with inertia disturbances.

\textbf{Theorem 1:} \textbf{[$\mathcal{H}_2$ norm for homogeneous power grids].} Consider an $N$-generator power system with both multiplicative and additive noise specified by parameters $\tilde{\Sigma} = (\tilde{A}_0, \tilde{A}_i, \tilde{B}, \tilde{C})$ that define the input-output mapping shown in \eqref{LTVstateform_reduced}. Suppose the inertia and damping are homogeneous, i.e., $\hat{M}_i = \hat{M}$, $\hat{\sigma}_i^2 = \hat{\sigma}^2, \forall i$ and $D = \beta I$. Consider also a general output matrix $\tilde{C} = \begin{bmatrix} \tilde{J}^{\frac{1}{2}} & 0\\
0 & \tilde{K}^{\frac{1}{2}}
\end{bmatrix}$, where matrices $\tilde{J}$ and $\tilde{K}$ are positive definite. Then the squared $\mathcal{H}_2$ norm is given by
\begin{equation}
\|\tilde{\Sigma}\|_{\mathcal{H}_2}^2 = \frac{1}{\hat{M}^2}\Tr\left[ \tilde{P}^{-1}\left(\hat{M}\tilde{J}\tilde{L}_B^{-1} + \tilde{K}\right) \right],
\end{equation}
where $\tilde{P} = \left[\left(\frac{2\beta}{\hat{M}}-\hat{\sigma}^2\beta^2\right)I-\hat{\sigma}^2\hat{M}\tilde{L}_B\right]$. 

\begin{proof}
See the Appendix.
\end{proof}

\textbf{Remark 2:  [\bm{$\mathcal{H}_2$} norm for specified outputs]}. We have the following expressions for the $\mathcal{H}_2$ norm of system \eqref{LTVstateform_reduced} for the three specific outputs mentioned in Section II:
\subsubsection{Phase cohesiveness} ($\tilde{C} = [\tilde{L}_{G}^{\frac{1}{2}}, 0]$)
\begin{equation*}
\|\tilde{\Sigma}\|_{\mathcal{H}_2}^2 = \frac{1}{\hat{M}}\Tr\left( \tilde{P}^{-1}\tilde{L}_G\tilde{L}_B^{-1} \right).
\end{equation*}
\subsubsection{Frequency} ($\tilde{C} = [0, I]$)
\begin{equation*}
\|\tilde{\Sigma}\|_{\mathcal{H}_2}^2 = \frac{1}{\hat{M}^2}\Tr\left( \tilde{P}^{-1}  \right).
\end{equation*}
\subsubsection{Phase cohesiveness \& frequency} ($\tilde{C} = \textrm{diag}\{\tilde{L}_G^{\frac{1}{2}},\kappa I\}$)

\begin{equation*}
\|\tilde{\Sigma}\|_{\mathcal{H}_2}^2 = \frac{1}{\hat{M}^2}\Tr\left[ \tilde{P}^{-1}\left(\hat{M}\tilde{L}_G\tilde{L}_B^{-1} + \kappa^2 I\right) \right].
\end{equation*}

\textbf{Corollary 1: [Second-moment bounded (Mean-square stability)  condition].} The power system in \eqref{LTVstateform_reduced} with both multiplicative and additive noise is second-moment bounded (mean-square stable) and has finite $\mathcal{H}_2$ norm if and only if
\begin{equation}\label{stabilitycondition}
\hat{\sigma}^2 < \frac{2\beta}{\hat{M}\left[\beta^2 + \lambda_{\textrm{max}}(\tilde{L}_B)\hat{M} \right]},
\end{equation}
where $\lambda_{\textrm{max}}(\tilde{L}_B)$ denotes the largest eigenvalue of Laplacian matrix $\tilde{L}_B$, and $\hat{\sigma}^2$ is the variance of the inverse distribution of $M$.

\begin{proof} Since $\tilde{A}_0$ is stable, the (standard) Lyapunov equation obtained when the inertia variance $\hat{\sigma}^2 = 0$ has a finite positive definite solution, and the corresponding $\mathcal{H}_2$ norm is finite. As $\hat{\sigma}^2$ increases, the smallest eigenvalue of the $\tilde{P}$ matrix defined in the $\mathcal{H}_2$ norm expression decreases (and thus largest eigenvalue of $\tilde{P}^{-1}$ increases), causing the $\mathcal{H}_2$ norm to increase. When $\hat{\sigma}^2$ approaches a critical value where $\tilde{P}$ goes from being positive definite to being singular, the $\mathcal{H}_2$ norm approaches infinity, and the system has unbounded second moment when the smallest eigenvalue of $\tilde{P}$ is zero. Examining the condition where the smallest eigenvalue of $\tilde{P}$ is zero yields the condition that guarantees bounded second moment.
\end{proof}

We will provide discussions of this Corollary in the context of our numerical results in the following section.

\section{Numerical results}

\begin{figure}[htbp!]
\centering
\includegraphics[width=2.5in]{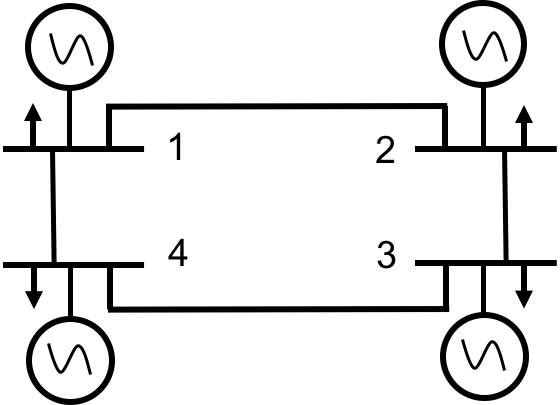}
\caption{A four-area interconnected power system.}
\label{fidg:fourbusnetwork}
\end{figure}

The results derived in the previous section indicate that the $\mathcal{H}_2$ norm of a power system with both multiplicative and additive noise depends on 1) the system topology (via the Laplacians $L_B$ and $L_G$); 2) the nominal system parameters (via the nominal inertia and damping coefficients    ); and 3) the distribution parameters of the multiplicative and additive disturbances. In this section, we present numerical simulations to analyze stability and to evaluate performance for three different outputs. Consider an interconnected power network with four areas (e.g., with homogeneous inertia $M$ and damping ratio $\beta$) shown in Fig.\ref{fidg:fourbusnetwork}. We assume this power network is Kron-reduced and single-phase equivalent, with line data given in Table \ref{table:LineParameters}.

\begin{table}[htbp!]
\centering
\caption{Line impedance parameters}
\begin{tabular}{|c|c|c|c|c|}
\hline
line $\mathcal{E}(i,j)$ & (1,2) & (2,3) & (3,4) & (4,1)\\
\hline
$r_{ij}$ & 0.4 & 0.5 & 0.6 & 0.28\\
\hline
$x_{ij}$ & 0.386 & 0.294 & 0.596 & 0.474\\
\hline
\end{tabular}
\label{table:LineParameters}
\end{table}

\begin{figure*}[htbp!]
\centering
\label{fidg:results}
\subfigure[phase cohesiveness output]{\label{fidg:results_phasecohesiveness} 
\includegraphics[width=2.25in]{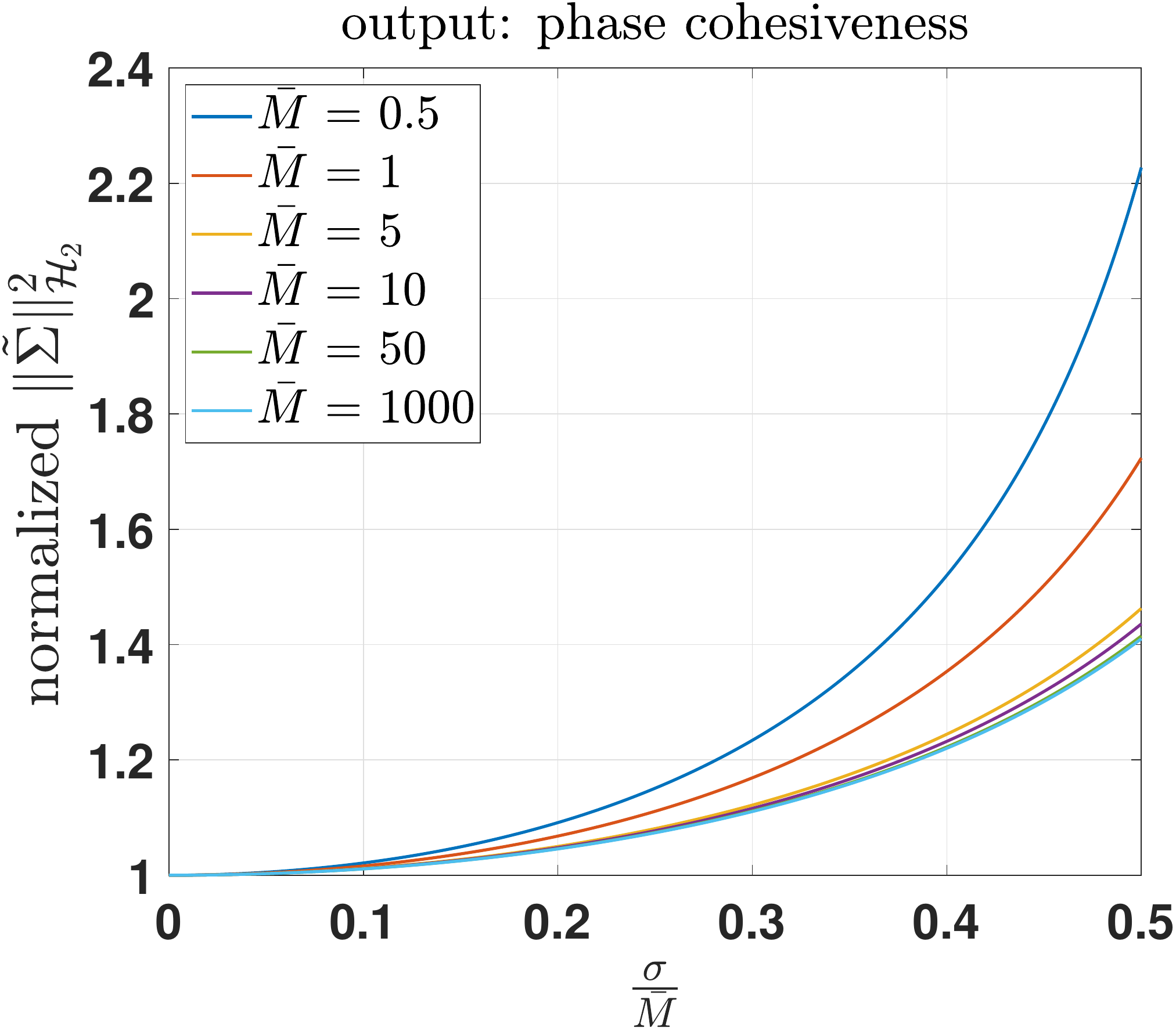}}
\subfigure[frequency output]{\label{fidg:results_frequency} 
\includegraphics[width=2.25in]{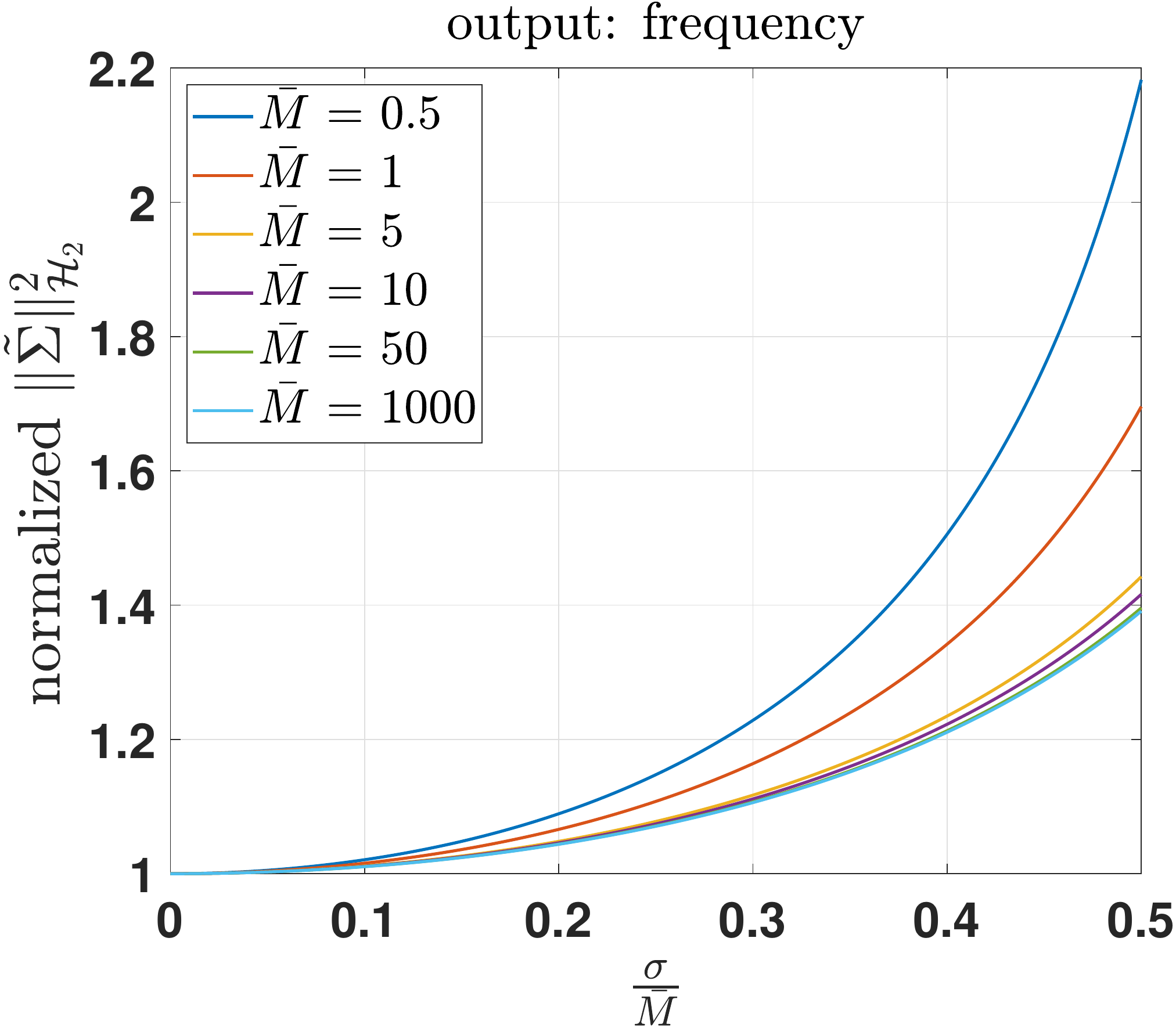}}
\subfigure[phase cohesiveness \& frequency output $\kappa = 10$]{\label{fidg:results_combined} 
\includegraphics[width=2.25in]{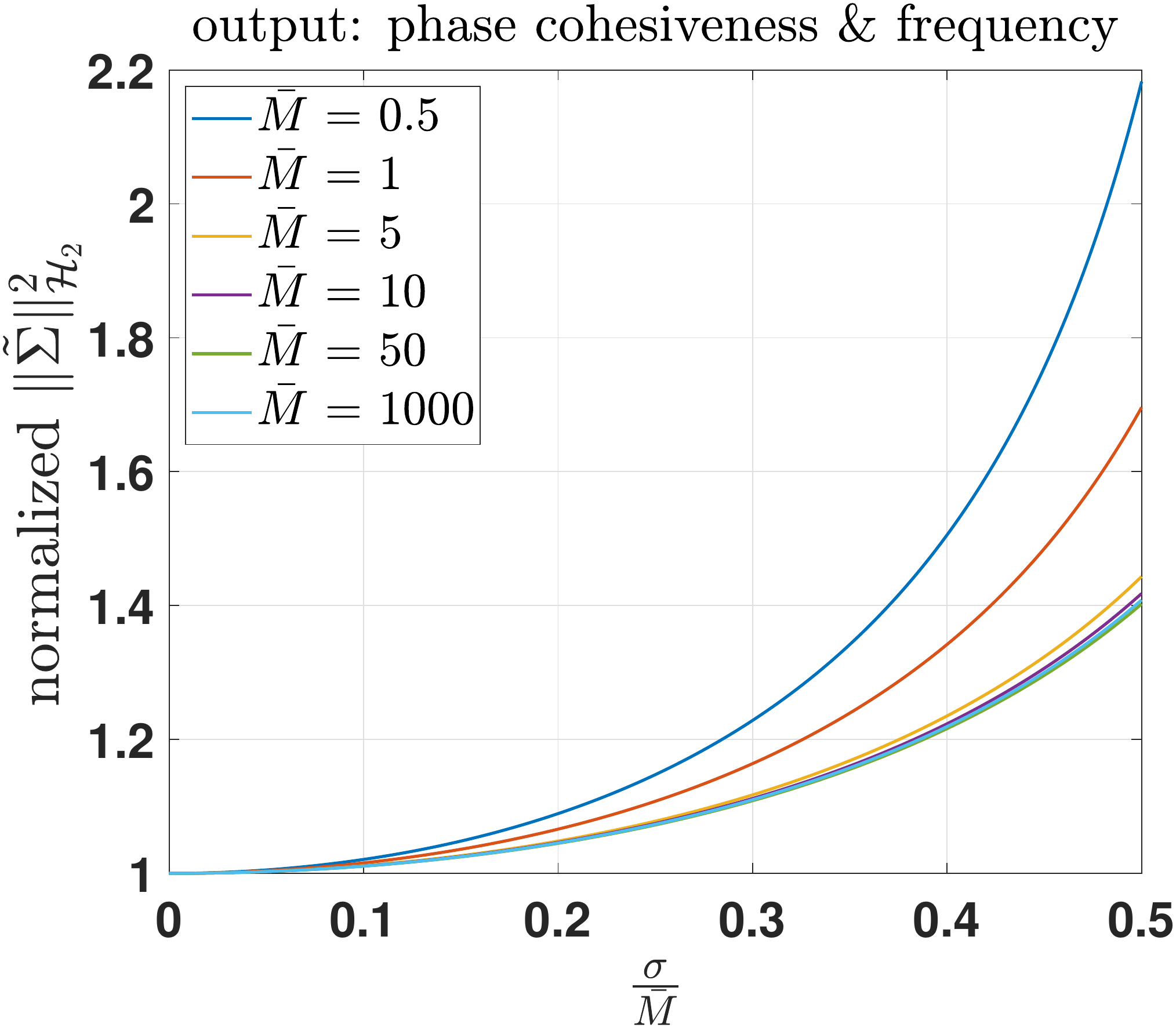}}
\caption{Comparison of the squared $\mathcal{H}_2$ norm of the stochastic system $\tilde{\Sigma}$ with three outputs under various values of $\bar{M}$ and inertia variance, quantified by $\frac{\sigma}{\bar{M}}$, for damping ratio $\beta = 1$. The results are normalized by the $\mathcal{H}_2$ norm with only additive noise, i.e., $\sigma^2 = 0$, [cf.  \cite{pirani2017system,tegling2015price}].
} 
\end{figure*}
\textbf{Remark 3: [Approximation of the inverse distribution of $\mathbf{M}$]}. To facilitate interpretation of insights derived from our previous analysis in terms of distribution parameters of $M$ instead of the inverse distribution $M^{-1}$, we seek to express the parameters (e.g., $\bar{M}$ and $\sigma^2$) of the original distribution $M$ in terms of parameters (e.g., $\hat{M}$ and $\hat{\sigma}^2$) of the inverse distribution $M^{-1}$. An approximate formula. commonly known as the $\delta$-method (based on a Taylor series expansion \cite{benaroya2005probability}), can be used to estimate the mean and variance of $M^{-1}$ considering variations around the mean value $\bar{M}$ \cite{mood1974introduction}. This yields an approximation of the distribution parameters of $M^{-1}$ give by
$\hat{M}^{-1} \approx \bar{M}^{-1} $, $\hat{\sigma}^2 \approx \frac{\sigma^2}{\bar{M}^4}$. The $\mathcal{H}_2$ norm of the system \eqref{LTVstateform_reduced} with the phase cohesiveness output is then approximated by
\begin{equation}\nonumber
\|\tilde{\Sigma}\|_{\mathcal{H}_2}^2 \approx \frac{1}{\bar{M}} \Tr\left( \left[ \left(\frac{2\beta}{\bar{M}}-\frac{\sigma^2\beta^2}{\bar{M}^4}\right)I - \frac{\sigma^2}{\bar{M}^3}\tilde{L}_B\right]^{-1}\tilde{L}_G\tilde{L}_B^{-1} \right).
\end{equation}
The $\mathcal{H}_2$ norm \eqref{LTVstateform_reduced} for the other outputs can be similarly approximated. We emphasize that this approximation is only to facilitate interpretation; it is possible but less intuitive to work with the inverse inertia distribution.

Fig.\ref{fidg:results_phasecohesiveness} shows the $\mathcal{H}_2$ norm for the phase cohesiveness output with increasing variance of inertia disturbances $\sigma^2$. In contrast to the additive noise case, where the $\mathcal{H}_2$ norm of \eqref{LTI ODE expression} is independent of the system inertia \cite{tegling2015price,pirani2017system}, the analytical (Theorem 1) and numerical results in Fig.\ref{fidg:results_phasecohesiveness} demonstrate that the $\mathcal{H}_2$ norm is strongly dependent on the nominal system inertia and its distribution. The results in Fig.\ref{fidg:results_phasecohesiveness} also demonstrate that a low-inertia power grid with larger inertia disturbance incur more resistive power losses in the presence of persistent disturbances or transient events.
Fig.\ref{fidg:results_frequency}-Fig.\ref{fidg:results_combined} shows that the  $\mathcal{H}_2$ norms of systems with frequency output and phase cohesiveness \& frequency output are also increasing functions of mean and variance of the system inertia disturbance. The system will suffer a larger stability degradation in the presence of stronger system inertia disturbance. The results from Fig.\ref{fidg:results_phasecohesiveness}-Fig.\ref{fidg:results_combined} also indicate that a power system with lower inertia is less robust to system inertia disturbances. In particular, the degradation of performance and stability margin is more severe in a power grid with less rotational inertia.

The analysis and numerical results indicate that the grid topology (via the Laplacians) plays an essential role in system stability. The first moment stability criteria (mean stability) requires the second smallest eigenvalue of Laplacian matrix $L_B$ to strictly larger than zero. Interestingly, the mean-square stability condition requires an upper bound on the \emph{largest} eigenvalue of Laplacian matrix $L_B$
\begin{equation}\nonumber
\lambda_{\textrm{max}}(L_B) < \frac{2\beta-\hat{M}\hat{\sigma}^2\beta^2}{\hat{\sigma}^2\hat{M}^2}.
\end{equation}
The mean and mean-square stability criteria provide a theoretical boundary on grid topology in terms of system coefficients and the distribution parameters of inertia disturbance (mean and variance). The mean-square stability condition indicates that highly connected power grids are more sensitive to multiplicative noise, which contrasts with the mean stability condition, where increasing algebraic connectivity improves stability robustness.

\section{Conclusion and outlooks}
In this paper, we proposed a stochastic swing equation with both multiplicative and additive noise to study low- and variable-inertia power system through the system $\mathcal{H}_2$ norm. The $\mathcal{H}_2$ norm can be computed by solving a generalized Lyapunov equation. For grids with homogeneous inertia and damping ratio, we derived an analytical expression of the $\mathcal{H}_2$ norm for various performance outputs. In contrast to the additive case, inertia variations may cause second-moment instability, even when the state mean is stable. Further, the performance metrics always depends on physical properties (via the nominal inertia, damping ratio, inertia distribution parameters) and network structure (via the susceptance matrix). Numerical results also indicate that the low-inertia grids are vulnerable to large system inertia disturbances. 

Ongoing works and potential future research directions include 1) more detailed stability analysis based on various performance metrics; 2) further numerical and analytical discussion of the stochastic swing equations and 3) 
design of optimal controllers for low- and variable-inertia grids with stochastic system inertia.



\section*{Appendix}
\begin{proof}
(of Theorem 1).
The squared $\mathcal{H}_2$ norm of \eqref{LTVstateform_reduced} is given by
\begin{equation}\label{H2tranceAppendix}
\|\tilde{\Sigma}\|_{\mathcal{H}_2}^2 = 
\Tr\left(\tilde{B}^\intercal \tilde{Q} \tilde{B}\right),
\end{equation}
where $\tilde{Q}$ is the observability Gramian, which can be obtained by solving the generalized Lyapunov equation $$\tilde{A}_0^\intercal\tilde{Q} + \tilde{Q}\tilde{A}_0 + \sum_{i=1}^N \hat{\sigma}_i^2\tilde{A}_i^\intercal \tilde{Q} \tilde{A}_i = -\tilde{C}^\intercal \tilde{C}.$$ The system matrices include the homogeneous inertia $\hat{M}$, inertia variance $\hat{\sigma}^2$, and damping coefficient $\beta$. We substitute $\tilde{B}$ in \eqref{H2tranceAppendix} and partition the Gramian as $\tilde{Q} = \begin{bmatrix}
\tilde{Q}_1 & \tilde{Q}_0\\
\tilde{Q}_0^\intercal & \tilde{Q}_2
\end{bmatrix}$. Due to the structure of the system matrices, the squared $\mathcal{H}_2$ norm of system $\tilde{\Sigma}$ becomes
\begin{equation}\label{H2normQ2}
\|\tilde{\Sigma}\|_{\mathcal{H}_2}^2 = \frac{1}{\hat{M}^2}\Tr\left( \tilde{Q}_2\right).
\end{equation}

Expanding the generalized Lyapunov equation yields
\begin{equation}
\scalebox{0.93}[1]{$
\begin{aligned}
\begin{bmatrix}
-\frac{1}{\hat{M}}\tilde{L}_B\tilde{Q}_0^\intercal, & -\frac{1}{\hat{M}}\tilde{L}_B\tilde{Q}_2\\
\tilde{Q}_1 - \frac{\beta}{\hat{M}}\tilde{Q}_0^\intercal, & \tilde{Q}_0-\frac{\beta}{\hat{M}}\tilde{Q}_2
\end{bmatrix}+
\begin{bmatrix}
-\frac{1}{\hat{M}}\tilde{Q}_0\tilde{L}_B, & \tilde{Q}_1 - \frac{\beta}{\hat{M}}\tilde{Q}_0\\
-\frac{1}{\hat{M}}\tilde{Q}_2\tilde{L}_B, & \tilde{Q}_0^\intercal-\frac{\beta}{\hat{M}}\tilde{Q}_2
\end{bmatrix} \\
+\hat{\sigma}^2
\begin{bmatrix}
\tilde{L}_B\tilde{Q}_2\tilde{L}_B, & \beta \tilde{L}_B\tilde{Q}_2\\
\beta \tilde{Q}_2\tilde{L}_B, & \beta^2 \tilde{Q}_2
\end{bmatrix}=
\begin{bmatrix}
-\tilde{J}, & 0\\
0, & -\tilde{K}
\end{bmatrix},
\end{aligned}
$}
\end{equation}
and the diagonal blocks are
\begin{subequations}
\begin{align}
-\frac{1}{\hat{M}}\tilde{L}_B\tilde{Q}_0^\intercal - \frac{1}{\hat{M}}\tilde{Q}_0\tilde{L}_B + \hat{\sigma}^2\tilde{L}_B\tilde{Q}_2\tilde{L}_B = -\tilde{J},\\
\tilde{Q}_0 - \frac{\beta}{\hat{M}}\tilde{Q}_2 + \tilde{Q}_0^\intercal - \frac{\beta}{\hat{M}}\tilde{Q}_2 + \hat{\sigma}^2\beta^2\tilde{Q}_2 = -\tilde{K}.
\end{align}
\end{subequations}
Since $\tilde{L}_B$ is nonsingular, the above two equations can be rearranged as
\begin{subequations}
\begin{align}
\label{LyapunovEq1}-\tilde{L}_B\tilde{Q}^\intercal_0\tilde{L}^{-1}_B-\tilde{Q}_0+\hat{\sigma}^2\hat{M}\tilde{L}_B\tilde{Q}_2=-\hat{M}\tilde{J}\tilde{L}_B^{-1} \\
\label{LyapunovEq2}\tilde{Q}_0+\tilde{Q}^\intercal_0+\left(\hat{\sigma}^2\beta^2-\frac{2\beta}{\hat{M}}\right)\tilde{Q}_2= - \tilde{K}.
\end{align}
\end{subequations}
Adding \eqref{LyapunovEq1} to \eqref{LyapunovEq2} and multiplying by $-1$ gives
\begin{equation}
-\tilde{Q}_0^\intercal+\tilde{L}_B\tilde{Q}^\intercal_0\tilde{L}_B^{-1} + \tilde{P}\tilde{Q}_2 = \hat{M}\tilde{J}\tilde{L}^{-1}_B + \tilde{K},
\end{equation}
where $\tilde{P} = \left[\left(\frac{2\beta}{\hat{M}}-\hat{\sigma}^2\beta^2\right)I-\hat{\sigma}^2\hat{M}\tilde{L}_B\right]$. Then multiplying by $\tilde{P}^{-1}$ and taking the trace gives
\begin{equation}\nonumber
\begin{aligned}
& \Tr\left(\tilde{P}^{-1}\left[ \tilde{L}_B\tilde{Q}_0^\intercal\tilde{L}_B^{-1} - \tilde{Q}_0^\intercal \right]\right)
+\Tr\left(\tilde{Q}_2\right)\\
& ~~~~~~~~~~~~~~~~~~~~~~~~~~~~~~~=\Tr\left[\tilde{P}^{-1}\left(\hat{M}\tilde{J}\tilde{L}_B^{-1} + \tilde{K}\right)\right].
\end{aligned}
\end{equation}
We will show that the first term in the above equation is zero. We define $\tilde{P} = \left[aI-b\tilde{L}_B\right]$, where $a = \frac{2\beta}{\hat{M}}-\hat{\sigma}^2\beta^2$ and $ b = \hat{\sigma}^2\hat{M}$. It can be seen that $\tilde{P}^{-1}$ and $\tilde{L}_B$ commute by expanding the term $\left[I - \frac{b}{a}\tilde{L}_B\right]^{-1}$ in a Neumann series
\begin{equation}
\begin{aligned}
\tilde{P}^{-1}\tilde{L}_B = \left[aI - b\tilde{L}_B\right]^{-1}\tilde{L}_B = \frac{1}{a}\left[I - \frac{b}{a}\tilde{L}_B\right]^{-1}\tilde{L}_B\\
= \frac{1}{a}\sum_{k=0}^{\infty}\left(\frac{b}{a}\tilde{L}_B\right)^k\tilde{L}_B = \tilde{L}_B\tilde{P}^{-1}.
\end{aligned}
\end{equation}
Since $\tilde{P}^{-1}\tilde{L}_B = \tilde{L}_B\tilde{P}^{-1}$, it follows that
\begin{equation}
\Tr\left(\tilde{P}^{-1}\left[ \tilde{L}_B\tilde{Q}_0^\intercal\tilde{L}_B^{-1} - \tilde{Q}_0^\intercal \right]\right) = 0,
\end{equation}
\begin{equation}
\Tr\left(\tilde{Q}_2\right)=\Tr\left[\tilde{P}^{-1}\left(\hat{M}\tilde{J}\tilde{L}_B^{-1} + \tilde{K}\right)\right].
\end{equation}
Finally, substituting $\Tr(\tilde{Q}_2)$ into \eqref{H2normQ2} leads to
\begin{equation}
\|\tilde{\Sigma}\|_{\mathcal{H}_2}^2 = \frac{1}{\hat{M}^2}\Tr\left[ \tilde{P}^{-1}\left(\hat{M}\tilde{J}\tilde{L}_B^{-1} + \tilde{K}\right) \right],
\end{equation}
which concludes the proof.
\end{proof}

\bibliography{reference}

\begin{thebibliography}{10}

\bibitem{guo2019performance}
Y.~Guo and T.~H. Summers, ``A performance and stability analysis of low-inertia
  power grids with stochastic system inertia,'' in {\em American Control
  Conference}, pp.~1--6, July 2019.

\bibitem{chiang1987foundations}
H.-D. Chiang, F.~Wu, and P.~Varaiya, ``Foundations of direct methods for power
  system transient stability analysis,'' {\em IEEE Transactions on Circuits and
  systems}, vol.~34, no.~2, pp.~160--173, 1987.

\bibitem{chiang2011direct}
H.-D. Chiang, {\em Direct methods for stability analysis of electric power
  systems: theoretical foundation, BCU methodologies, and applications}.
\newblock John Wiley \& Sons, 2011.

\bibitem{kundur2004definition}
P.~Kundur {\em et~al.}, ``Definition and classification of power system
  stability,'' {\em IEEE Transactions on Power Systems}, vol.~19, no.~2,
  pp.~1387--1401, 2004.

\bibitem{bergen1981structure}
A.~R. Bergen and D.~J. Hill, ``A structure preserving model for power system
  stability analysis,'' {\em IEEE Transactions on Power Apparatus and Systems},
  no.~1, pp.~25--35, 1981.

\bibitem{tegling2015price}
E.~Tegling, B.~Bamieh, and D.~F. Gayme, ``The price of synchrony: Evaluating
  the resistive losses in synchronizing power networks,'' {\em IEEE
  Transactions Control of Network Systems}, vol.~2, no.~3, pp.~254--266, 2015.

\bibitem{pirani2017system}
M.~Pirani, J.~W. Simpson-Porco, and B.~Fidan, ``System-theoretic performance
  metrics for low-inertia stability of power networks,'' in {\em 56th IEEE
  Annual Conference on Decision and Control}, pp.~5106--5111, 2017.

\bibitem{motter2013spontaneous}
A.~E. Motter, S.~A. Myers, M.~Anghel, and T.~Nishikawa, ``Spontaneous synchrony
  in power-grid networks,'' {\em Nature Physics}, vol.~9, no.~3, p.~191, 2013.

\bibitem{dorfler2012synchronization}
F.~D{\"o}rfler and F.~Bullo, ``Synchronization and transient stability in power
  networks and nonuniform kuramoto oscillators,'' {\em SIAM Journal on Control
  and Optimization}, vol.~50, no.~3, pp.~1616--1642, 2012.

\bibitem{luo2018coordinated}
H.~Luo, Z.~Hu, H.~Zhang, and H.~Chen, ``Coordinated active power control
  strategy for deloaded wind turbines to improve regulation performance in
  {AGC},'' {\em IEEE Transactions on Power Systems}, 2018.

\bibitem{ahmadyar2018framework}
A.~S. Ahmadyar, S.~Riaz, G.~Verbi{\v{c}}, A.~Chapman, and D.~J. Hill, ``A
  framework for assessing renewable integration limits with respect to
  frequency performance,'' {\em IEEE Transactions on Power Systems}, vol.~33,
  no.~4, pp.~4444--4453, 2018.

\bibitem{ulbig2013impact}
A.~Ulbig, T.~S. Borsche, and G.~Andersson, ``Impact of low rotational inertia
  on power system stability and operation,'' {\em IFAC Proceedings Volumes},
  vol.~47, no.~3, pp.~7290--7297, 2014.

\bibitem{ulbig2015analyzing}
A.~Ulbig, T.~S. Borsche, and G.~Andersson, ``Analyzing rotational inertia, grid
  topology and their role for power system stability,'' {\em
  IFAC-PapersOnLine}, vol.~48, no.~30, pp.~541--547, 2015.

\bibitem{kroposki2017achieving}
B.~Kroposki, B.~Johnson, Y.~Zhang, V.~Gevorgian, P.~Denholm, B.-M. Hodge, and
  B.~Hannegan, ``Achieving a 100\% renewable grid: Operating electric power
  systems with extremely high levels of variable renewable energy,'' {\em IEEE
  Power and Energy Magazine}, vol.~15, no.~2, pp.~61--73, 2017.

\bibitem{winter2015pushing}
W.~Winter, K.~Elkington, G.~Bareux, and J.~Kostevc, ``Pushing the limits:
  Europe's new grid: Innovative tools to combat transmission bottlenecks and
  reduced inertia,'' {\em IEEE Power and Energy Magazine}, vol.~13, no.~1,
  pp.~60--74, 2015.

\bibitem{tielens2012grid}
P.~Tielens and D.~Van~Hertem, ``Grid inertia and frequency control in power
  systems with high penetration of renewables,'' 2012.

\bibitem{beck2007virtual}
H.-P. Beck and R.~Hesse, ``Virtual synchronous machine,'' in {\em 9th IEEE
  International Conference on Electrical Power Quality and Utilisation},
  pp.~1--6, 2007.

\bibitem{d2013virtual}
S.~D'Arco and J.~A. Suul, ``Virtual synchronous machines classification of
  implementations and analysis of equivalence to droop controllers for
  microgrids,'' in {\em IEEE Grenoble PowerTech}, pp.~1--7, 2013.

\bibitem{zhong2012control}
Q.-C. Zhong and T.~Hornik, {\em Control of power inverters in renewable energy
  and smart grid integration}, vol.~97.
\newblock John Wiley \& Sons, 2012.

\bibitem{dorfler2016breaking}
F.~D{\"o}rfler, J.~W. Simpson-Porco, and F.~Bullo, ``Breaking the hierarchy:
  Distributed control and economic optimality in microgrids,'' {\em IEEE
  Transactions on Control of Network Systems}, vol.~3, no.~3, pp.~241--253,
  2016.

\bibitem{bevrani2014virtual}
H.~Bevrani, T.~Ise, and Y.~Miura, ``Virtual synchronous generators: A survey
  and new perspectives,'' {\em International Journal of Electrical Power \&
  Energy Systems}, vol.~54, pp.~244--254, 2014.

\bibitem{zhu2018stability}
L.~Zhu and D.~J. Hill, ``Stability analysis of power systems: A network
  synchronization perspective,'' {\em SIAM Journal on Control and
  Optimization}, vol.~56, no.~3, pp.~1640--1664, 2018.

\bibitem{song2017network}
Y.~Song, D.~J. Hill, and T.~Liu, ``Network-based analysis of small-disturbance
  angle stability of power systems,'' {\em IEEE Transactions on Control of
  Network Systems}, 2017.

\bibitem{colombino2017global}
M.~Colombino, D.~Gro{\ss}, J.-S. Brouillon, and F.~D{\"o}rfler, ``Global phase
  and magnitude synchronization of coupled oscillators with application to the
  control of grid-forming power inverters,'' {\em IEEE Transactions on
  Automatic Control}, 2019.

\bibitem{torres2015synchronization}
L.~A. T{\^o}rres, J.~P. Hespanha, and J.~Moehlis, ``Synchronization of
  identical oscillators coupled through a symmetric network with dynamics: A
  constructive approach with applications to parallel operation of inverters,''
  {\em IEEE Transactions on Automatic Control}, vol.~60, no.~12,
  pp.~3226--3241, 2015.

\bibitem{sinha2017uncovering}
M.~Sinha, F.~D{\"o}rfler, B.~B. Johnson, and S.~V. Dhople, ``Uncovering droop
  control laws embedded within the nonlinear dynamics of van der pol
  oscillators,'' {\em IEEE Transactions on Control of Network Systems}, vol.~4,
  no.~2, pp.~347--358, 2017.

\bibitem{zhao2014design}
C.~Zhao, U.~Topcu, N.~Li, and S.~Low, ``Design and stability of load-side
  primary frequency control in power systems,'' {\em IEEE Transactions on
  Automatic Control}, vol.~59, no.~5, pp.~1177--1189, 2014.

\bibitem{mallada2017optimal}
E.~Mallada, C.~Zhao, and S.~Low, ``Optimal load-side control for frequency
  regulation in smart grids,'' {\em IEEE Transactions on Automatic Control},
  vol.~62, no.~12, pp.~6294--6309, 2017.

\bibitem{curi2017control}
S.~Curi, D.~Gro{\ss}, and F.~D{\"o}rfler, ``Control of low-inertia power grids:
  A model reduction approach,'' in {\em 56th IEEE Annual Conference on Decision
  and Control}, pp.~5708--5713, 2017.

\bibitem{li2016connecting}
N.~Li, C.~Zhao, and L.~Chen, ``Connecting automatic generation control and
  economic dispatch from an optimization view,'' {\em IEEE Transactions on
  Control of Network Systems}, vol.~3, no.~3, pp.~254--264, 2016.

\bibitem{mevsanovic2016comparison}
A.~Me{\v{s}}anovi{\'c}, U.~M{\"u}nz, and C.~Heyde, ``Comparison of
  $\mathcal{H}_\infty$, $\mathcal{H}_2$, and pole optimization for power system
  oscillation damping with remote renewable generation,'' {\em
  IFAC-PapersOnLine}, vol.~49, no.~27, pp.~103--108, 2016.

\bibitem{ahmadyar2016assessment}
A.~S. Ahmadyar, S.~Riaz, G.~Verbi{\v{c}}, J.~Riesz, and A.~Chapman,
  ``Assessment of minimum inertia requirement for system frequency stability,''
  in {\em IEEE International Conference on Power System Technology}, pp.~1--6,
  Sept. 2016.

\bibitem{poolla2017optimal}
B.~K. Poolla, S.~Bolognani, and F.~D{\"o}rfler, ``Optimal placement of virtual
  inertia in power grids,'' {\em IEEE Transactions on Automatic Control},
  vol.~62, no.~12, pp.~6209--6220, 2017.

\bibitem{dorfler2013kron}
F.~D{\"o}rfler and F.~Bullo, ``Kron reduction of graphs with applications to
  electrical networks.,'' {\em IEEE Transactions on Circuits and Systems},
  vol.~60, no.~1, pp.~150--163, 2013.

\bibitem{purchala2005usefulness}
K.~Purchala, L.~Meeus, D.~Van~Dommelen, and R.~Belmans, ``Usefulness of {DC}
  power flow for active power flow analysis,'' in {\em IEEE PES General
  Meeting}, pp.~454--459, 2005.

\bibitem{oksendal2003stochastic}
B.~{\O}ksendal, ``Stochastic differential equations,'' in {\em Stochastic
  differential equations}, pp.~65--84, Springer, 2003.

\bibitem{jiang2017performance}
Y.~Jiang, R.~Pates, and E.~Mallada, ``Performance tradeoffs of dynamically
  controlled grid-connected inverters in low inertia power systems,'' in {\em
  56th IEEE Annual Conference on Decision and Control}, pp.~5098--5105, 2017.

\bibitem{damm2004rational}
T.~Damm, {\em Rational matrix equations in stochastic control}, vol.~297.
\newblock Springer Science \& Business Media, 2004.

\bibitem{verriest2008time}
E.~Verriest, ``Time variant balancing and nonlinear balanced realizations,'' in
  {\em Model Order Reduction: Theory, Research Aspects and Applications},
  pp.~213--250, Springer, 2008.

\bibitem{zhang2002h2}
L.~Zhang and J.~Lam, ``On h2 model reduction of bilinear systems,'' {\em
  Automatica}, vol.~38, no.~2, pp.~205--216, 2002.

\bibitem{kleinman1969stability}
D.~Kleinman, ``On the stability of linear stochastic systems,'' {\em IEEE
  Transactions on Automatic Control}, vol.~14, no.~4, pp.~429--430, 1969.

\bibitem{benner2011lyapunov}
P.~Benner and T.~Damm, ``Lyapunov equations, energy functionals, and model
  order reduction of bilinear and stochastic systems,'' {\em SIAM Journal on
  Control and Optimization}, vol.~49, no.~2, pp.~686--711, 2011.

\bibitem{damm2008direct}
T.~Damm, ``Direct methods and adi-preconditioned krylov subspace methods for
  generalized lyapunov equations,'' {\em Numerical Linear Algebra with
  Applications}, vol.~15, no.~9, pp.~853--871, 2008.

\bibitem{bai2006projection}
Z.~Bai and D.~Skoogh, ``A projection method for model reduction of bilinear
  dynamical systems,'' {\em Linear algebra and its applications}, vol.~415,
  no.~2-3, pp.~406--425, 2006.

\bibitem{boyd1994linear}
S.~Boyd, L.~El~Ghaoui, E.~Feron, and V.~Balakrishnan, {\em Linear matrix
  inequalities in system and control theory}, vol.~15.
\newblock SIAM, 1994.

\bibitem{benaroya2005probability}
H.~Benaroya, S.~M. Han, S.~Han, and M.~Nagurka, {\em Probability models in
  engineering and science}.
\newblock CRC press, 2005.

\bibitem{mood1974introduction}
A.~Mood, F.~Graybill, and D.~Boes, ``Introduction to statistical theory,''
  1974.

\end{thebibliography}
\bibliographystyle{ieeetr}

\end{document}